*****LATEX FILE*****

\documentstyle[12pt]{article}

\def\fr{\frac}

\def\Si{\Sigma}

\def\de{\delta}

\def\a{a^{\dagger}}
\def\e{\exp_q}

\newcommand{\bd}{\begin{displaymath}}
\newcommand{\ed}{\end{displaymath}}
\newcommand{\bb}{\begin{equation}}
\newcommand{\ee}{\end{equation}}

\begin{document}
\baselineskip 1.8 \baselineskip


\vspace{.2cm}

\begin{center}
\Large {\bf $gl_q(n)$-covariant Oscillators and q-Deformed Quantum Mechanics
in n Dimensions}
\\[1cm]
\large W.-S.Chung \\[.3cm]
\normalsize  
Theory Group, Department of Physics, College of Natural Sciences,  \\
\normalsize  Gyeongsang National University,   \\
\normalsize   Jinju, 660-701, Korea
\end{center}

\vspace{0.5cm}
\begin{abstract}
In this paper the coherent state for $gl_q(n)$-covariant oscillators is
constructed and
is shown to satisfy the completeness relation. And the q-analogue of quantum
mechanics
in n dimensions is obtained by using $gl_q(n)$ oscillators.
\end{abstract}

\setcounter{page}{1}
\section{Introduction}

Quantum                                                               
groups                                                                      
or                        
q-deformed                    
Lie         
algebra                
implies   
some    
specific 
deformations 
of                    classical                      Lie           algebras.

From           
a             
mathematical      
point   
of  
view,     
it 
is  
a 
non-commutative                                                  
associative                                                                 
Hopf         
algebra.         
The       
structure        
and 
representation 
theory    of                                   quantum               groups  
have                 
been             
developed        
extensively  
by  
Jimbo    
[1] 
and 
Drinfeld                                                                
[2].                                                                        
            
The                                   
q-deformation                    
of        
Heisenberg       
algebra   
was     
made  
by 
Arik and Coon [3], Macfarlane [4] and Biedenharn [5].
Recently                                                               
there                                               
has                           
been                    
some              
interest               
in    
more     
general   
deformations 
involving                                                                 
an                                                                          
arbitrary                      
real                           
functions     
of           
weight    
generators  
and   
including 
q-deformed algebras as a special case [6-10].

In the mean time some theoretical physicists studied the q-deformation of
quantum mechanic in one dimension [11-16].
The purpose of this paper is to use the $gl_q(n)$-covariant oscillator
algebra to
construct the q-analogue of the quantum mechanics with harmonic potential
in n dimensions.
\section{Coherent states of $gl_q(n)$-covariant oscillator algebra}

$gl_q(n)$-covariant oscillator algebra is defined as [17]
\bd
\a_i\a_j =\sqrt{q} \a_j \a_i,~~~(i<j)  
\ed
\bd
a_ia_j=\fr{1}{\sqrt{q}}a_j a_i,~~~(i<j)
\ed
\bd
a_i\a_j=\sqrt{q} \a_ja_i,~~~(i \neq j)
\ed
\bd
a_i\a_i =1+q \a_ia_i +(q-1) \Si_{k=i+1}^n\a_k a_k,~~~(i=1,2,\cdots,n-1)
\ed
\bd
a_n \a_n =1+q \a_n a_n,
\ed
\bb
[N_i, a_j]=-\de_{ij}a_j,~~~[N_i, \a_j]=\de_{ij}\a_j,~~~(i,j=1,2,\cdots, n )
\ee
where we restrict our concern to the case that $q$ is real and $0<q<1$.
Here $N_i$ plays a role of number operator and $a_i(\a_i)$ plays a role of
annihilation(creation) operator.
From the above algebra one can obtain the relation between the number operators
and mode opeartors as follows
\bb
\a_ia_i=q^{\Si_{k=i+1}^nN_k}[N_i],
\ee
where $[x]$ is called a q-number and is defined as
\bd
[x]=\fr{q^{x}-1}{q-1}.
\ed
\def\nn{|n_1,n_2,\cdots,n_n>}
Let us introduce the Fock space basis $\nn$ for the number operators
$N_1,N_2,\cdots, N_n$ satisfying
\bb
N_i\nn=n_i\nn,~~~(n_1,n_2,\cdots,n_n=0,1,2\cdots)
\ee
Then we have the following representation
\bd
a_i\nn=\sqrt{q^{\Si_{k=i+1}^nn_k}[n_i]}|n_1,\cdots, n_i-1,\cdots,n_n>
\ed
\bb
\a_i\nn=\sqrt{q^{\Si_{k=i+1}^nn_k}[n_i+1]}|n_1,\cdots, n_i+1,\cdots,n_n>.
\ee

From the above representation we know that there exists the ground state
$|0,0,\cdots,0>$ satisfying
$a_i|0,0>=0$ for all $i=1,2,\cdots,n$. Thus the state $\nn$ is obtatind by
applying the creation operators
to the ground state $|0,0,\cdots,0>$
\bb
\nn=\fr{(\a_n)^{n_n}\cdots(\a_1)^{n_1}}{\sqrt{[n_1]!\cdots
[n_n]!}}|0,0,\cdots,0>.
\ee
If we introduce the scale operators as follows
\bb
Q_i=q^{N_i},~~(i=1,2,\cdots,n),
\ee
we have from the algebra (1)
\bb
[a_i,\a_i]=Q_iQ_{i+1}\cdots Q_n.
\ee
Acting the operators $Q_i$'s on the basis $\nn$ produces
\bb
Q_i\nn=q^{n_i}\nn .
\ee

From the relation $a_i a_j =\fr{1}{\sqrt{q}}a_j a_i,~~(i<j)$, the coherent
states for $gl_q(n)$ algebra
is defined as
\bb
a_i|z_1,\cdots,z_i,\cdots,z_n>=z_i|z_1,\cdots,
z_{i},\sqrt{q}z_{i+1},\cdots,\sqrt{q}z_n>.
\ee
Solving the eq.(9) we obtain
\bb
|z_1,z_2,\cdots,z_n>=c(z_1,\cdots,z_n)\Si_{n_1,n_2,\cdots,n_n=0}^{\infty} 
\fr{z_1^{n_1}z_2^{n_2}\cdots z_n^{n_n}}{\sqrt{[n_1]![n_2]!\cdots [n_n]!}}\nn .
\ee
Using eq.(5) we can rewrite eq.(10) as
\bb
|z_1,z_2,\cdots,z_n>=c(z_1,\cdots,z_n)
\e(z_n\a_n)\cdots\e(z_2\a_2)\e(z_1\a_1)|0,0,\cdots,0>.
\ee
where q-exponential function is defined as
\bd
\e(x)=\Si_{n=0}^{\infty}\fr{x^n}{[n]!}.
\ed
The q-exponential function satisfies the following recurrence relation
\bb
\e(q x)=[1-(1-q)x]\e(x)
\ee
Using the above relation and the fact that $0<q<1$, we obtain the formula
\bb
\e(x) =\Pi_{n=0}^{\infty}\fr{1}{1-(1-q)q^{n}x}
\ee
Using the normalization of the coherent state , we have
\bb
c(z_1,z_2,\cdots,z_n)=\e(|z_1|^2)\e(|z_2|^2)\cdots \e(|z_n|^2).
\ee
The coherent state satisfies the completeness relation
\bb
\int\cdots \int
|z_1,z_2,\cdots,z_n><z_1,z_2,\cdots,z_n|\mu(z_1,z_2,\cdots,z_n) d^2z_1
d^2z_2\cdots d^2 z_n=I,
\ee
where the weighting function $\mu(z_1,z_2,\cdots,z_n)$ is defined as
\bb
\mu(z_1,z_2,\cdots,z_n)=\fr{1}{\pi^2}\Pi_{i=1}^n\fr{\e(|z_i|^2)}
{\e(q|z_i|^2)}.
\ee
In deriving eq.(15) we used the formula
\bb
\int_0^{1/(1-q)}x^n \e(q x)^{-1} d_q x=[n]!
\ee

\section{q-Deformed Weyl-Heisenberg Group}
The purpose of this section is to explain what is the q-analogue of the
q-deformed
Weyl-Heisenberg group.
From the algebra (1) we obtain
\bd
a_i f(\a_i)=f(q \a_i) a_i +(Df)(\a_i) Q_{i+1}\cdots Q_n
\ed
\bb
a_n f(\a_n)=f(q \a_n) a_n +(Df)(\a_n) ,
\ee
where $D$ is called q-derivative and defined as
\bd
DF(x)=\fr{F(x)-F(qx)}{x(1-q)}.
\ed
Putting $f(x)=\e(tx)$ we have
\bb
a_i \e(t\a_i) =\e(qt\a_i)a_i +t\e(t\a_i)Q_{i+1}\cdots Q_n.
\ee

Using the formula (12) we have
\bb
a_i^n\e(t\a_i) =\e(t\a_i)(a_i+tQ_iQ_{i+1}\cdots Q_n)^n
\ee
and thus
\bb
\e(s_ia_i)\e(t_i\a_i)=\e(t_i\a_i)\e(s_ia_i+s_it_iQ_iQ_{i+1}\cdots Q_n).
\ee
Taking account of $[a_i,Q_i]_q=a_iQ_i -q Q_i a_i=0$,we have 
\bb
\e(s_ia_i)\e(t_i\a_i)=\e(t_i\a_i)\e(s_it_iQ_iQ_{i+1}\cdots Q_n)\e(s_ia_i).
\ee
If we muliply above equations from $i=1$ to $n$, we obtain the q-deformed
Weyl-Heisenberg
 relation
\bb
\Pi_{i=1}^n\e(s_ia_i)\e(t_i\a_i)=
\Pi_{i=1}^n\e(t_i\a_i)\e(s_it_iQ_iQ_{i+1}\cdots Q_n)\e(s_ia_i).
\ee

\section{q-deformed quantum mechanics in n dimensions}
It is intersting to study the q-deformed harmonic oscillator system in n
dimensions.
In order to formulate it we define the position and momentum operators 
\bd
X_i=\fr{1}{\sqrt{2}}(a_i +\a_i)
\ed
\bb
P_i=-\fr{i}{\sqrt{2}}(a_i -\a_i).
\ee
Then the Hamiltonian of q-deformed harmonic oscillator in n dimensions is
given by
\bb
H=\Si_{i=1}^n H_i,
\ee
where
\bb
H_i=\fr{1}{2}(P_i^2 +X_i^2) =\fr{1}{2} (a_i\a_i+\a_i a_i).
\ee
Now, the q-cannonical commutation relation can be expressed by
\bb
X_iP_i-P_iX_i =i(\fr{q+1}{2})^{i-n-1}
+i(q-1)\Si_{k=i}^n (\fr{q+1}{2})^{i-k-1}H_k.
\ee
Expressing $H_i$'s in terms of $Q_i$'s operators, we get
\bd
H_i=\fr{q+1}{2(q-1)}Q_iQ_{i+1}\cdots Q_n -\fr{1}{q-1}Q_{i+1}Q_{i+2}\cdots Q_n,
~~(i=1,2,\cdots, n-1)
\ed
\bb
H_n=\fr{q+1}{2(q-1)}Q_n -\fr{1}{q-1}.
\ee
Thus the Hamiltonian is given by
\bb
H=\fr{Q-1}{q-1}+\fr{1}{2}\Si_{i=1}^nQ_iQ_{i+1}\cdots Q_n
\ee
where
\bd
Q=Q_1Q_2 \cdots Q_n
\ed
Thus we have
\bb
H|n_1,\cdots,n_n>=E(n_1,\cdots,n_n)|n_1,\cdots,n_n>
\ee
where the energy spectrum is given by
\bb
E(n_1,\cdots,n_2)=[n_1+\cdots+n_n]+\fr{1}{2}\Si_{i=1}^nq^{n_1+\cdots +n_n}
\ee
\section{Concluding Remark}
In this paper we used $gl_q(n)$-covariant oscillator algebra to obtain its
coherent
state and showed  the completeness relation. Moreover we construct the
q-deformed
quantum mechanical hamiltonian in n dimensions by using $gl_q(n)$-covariant
oscillators.
In conclusion, it was known that we can obtain the q-analogue of
n-dimensional Schroedinger equation with harmonic potential
by using $gl_q(n)$-covariant oscillator system.

\section*{Acknowledgement}
This                   paper                was
supported         by  
the   KOSEF (961-0201-004-2)   
and   the   present   studies    were   supported   by   Basic  
Science 
Research Program, Ministry of Education, 1995 (BSRI-95-2413).

\end{document}